\documentclass[reqno,12pt]{amsart}
\usepackage{geometry}
\usepackage{amsmath,amstext, amsthm, empheq, extpfeil,  
amsbsy,amssymb,marvosym,fancyhdr,graphicx,amscd,amsfonts,latexsym,delarray,stackrel,lineno,color,cite,appendix,xcolor,url,hyperref}
\usepackage{wasysym}
\usepackage{subfig}

\usepackage{srcltx}
\usepackage{mathrsfs}
\usepackage{float} 
\usepackage[all]{xy}
\usepackage{t1enc}
\usepackage{mathrsfs}
\usepackage{pifont}

\usepackage{mathpple}
\usepackage[T1]{fontenc}

\definecolor{Green}{rgb}{0,1,0}
\definecolor{Blue}{RGB}{0,0,191}
\definecolor{mathmodecolor}{RGB}{0,102,0}
\definecolor{keywordcolor}{RGB}{0,51,151}
\definecolor{sourcebackgroundcolor}{RGB}{255,247,223}
\definecolor{unixagred}{RGB}{255,0,0}
\definecolor{lightgray}{RGB}{191,191,191}
\definecolor{green}{RGB}{1,191,191}

\newcommand*\patchAmsMathEnvironmentForLineno[1]{%
  \expandafter\let\csname old#1\expandafter\endcsname\csname #1\endcsname
  \expandafter\let\csname oldend#1\expandafter\endcsname\csname end#1\endcsname
  \renewenvironment{#1}%
     {\linenomath\csname old#1\endcsname}%
     {\csname oldend#1\endcsname\endlinenomath}}%
\newcommand*\patchBothAmsMathEnvironmentsForLineno[1]{%
  \patchAmsMathEnvironmentForLineno{#1}%
  \patchAmsMathEnvironmentForLineno{#1*}}%
\AtBeginDocument{%
\patchBothAmsMathEnvironmentsForLineno{equation}%
\patchBothAmsMathEnvironmentsForLineno{align}%
\patchBothAmsMathEnvironmentsForLineno{flalign}%
\patchBothAmsMathEnvironmentsForLineno{alignat}%
\patchBothAmsMathEnvironmentsForLineno{gather}%
\patchBothAmsMathEnvironmentsForLineno{multline}%
}
\newtheorem{thm}{Theorem}[section]

\newtheorem{lem}[thm]{Lemma}

\def\Sp{{\rm Spec}\,}

\def\Tr{{\rm Tr}}

\def\A{{\mathbb A}}

\def\Q{{\mathbb Q}}
\def\R{{\mathbb R}}

\def\Tr{{\rm Tr}}

\def\cS{{\mathcal S}}

\newcommand{\ie}{{\it i.e.\/}\ }

 \def\scal2{{\mathscr S}}

\def\sr0{{\cS^{\rm ev}_0}}

\def\sar0{{\cS_0(\A_\Q)}}

\newcommand{\nil}[1]{}

\author{Alain Connes}

 \title{\LARGE Heat expansion and zeta}
 \address{A.~Connes: Coll\`ege de France \\
3, rue d'Ulm \\ Paris, F-75005 France\\
I.H.E.S. 35 Route de Chartres Bures sur Yvette, F-91440} \email{alain\@@connes.org}


\begin{document}


\maketitle	

\begin{abstract}
We compute the full asymptotic expansion of the heat kernel $	\Tr(\exp(-tD^2))$ where $D$ is, assuming RH, the self-adjoint operator whose spectrum is formed of the imaginary parts of non-trivial zeros of the Riemann zeta function. The coefficients of the expansion are explicit expressions involving Bernoulli and Euler numbers. We relate the divergent terms with the heat kernel expansion of the Dirac square root of the prolate wave operator investigated in our joint work with Henri Moscovici.
\end{abstract}
\vspace{1cm}
\begin{center} % Centers your dedication text
\itshape % Italicizes the dedication text
Dedicated to Fedor Sukochev %\\ for inspiring me to reach further.
\end{center}
\tableofcontents
\section{Introduction}
We compute in closed form, using Bernoulli and Euler numbers, the  full asymptotic expansion of the heat trace $	\Tr(\exp(-tD^2))$ where $D$ is, assuming RH, the self-adjoint operator whose spectrum is formed of the imaginary parts of non-trivial zeros of the Riemann zeta function.
The heat expansion, \ie the asymptotic expansion for $t\to 0$, of the trace
	$	\Tr(\exp(-tD^2))$ where $D$ is a  self-adjoint operator, plays a key role in the theory of spectral asymptotics. 
It is a powerful tool in mathematical physics, differential geometry, and spectral theory, describing how solutions to the heat equation on a manifold or more general spaces behave, especially in the context of the Laplace operator or more general elliptic operators. The expansion typically takes the form of an asymptotic series, with coefficients that reflect geometric and topological properties of the space. These coefficients are often referred to as Seeley-DeWitt or heat kernel coefficients. This type of asymptotic expansion also plays a key role in noncommutative geometry as a way to extend curvature invariants to noncommutative spaces such as noncommutative tori, and Fedor Sukochev has greatly contributed to this research together with his collaborators with whom he forms a very active research group.

Computing the coefficients of the asymptotic expansion in closed form is a challenging problem, and there are only a few well-understood cases where all the coefficients can be explicitly calculated. Some of these examples include:

1. Flat Spaces: In Euclidean spaces $\left(\mathbb{R}^n\right)$ or flat tori, the heat kernel and its coefficients can be computed exactly because the geometry is simple and fully understood. The heat kernel in these cases is given by explicit formulas involving exponential functions.

2. Spheres: For the standard sphere $\left(S^n\right.$ ) in various dimensions, the heat kernel coefficients can be computed due to the high symmetry of the space and the explicit knowledge of the eigenvalues and eigenfunctions of the Laplacian.

3. Compact Lie Groups: For compact Lie groups equipped with a biinvariant metric, the heat kernel can be expressed in terms of the sum over the group's representations. This allows for the explicit calculation the heat kernel coefficients due to the algebraic structure of the group and the representation theory.

4. Certain Symmetric Spaces: Symmetric spaces of compact type can sometimes allow for the explicit computation of heat kernel coefficients, again due to the symmetry and the structure of the space allowing for the representation of the heat kernel in a manageable form.

5. Hyperbolic Spaces: For constant negative curvature spaces  (hyperbolic spaces), the heat kernel can be expressed in terms  of special functions, and in some cases, the coefficients can be computed explicitly.

\begin{figure}[H]	\begin{center}
\includegraphics[scale=1]{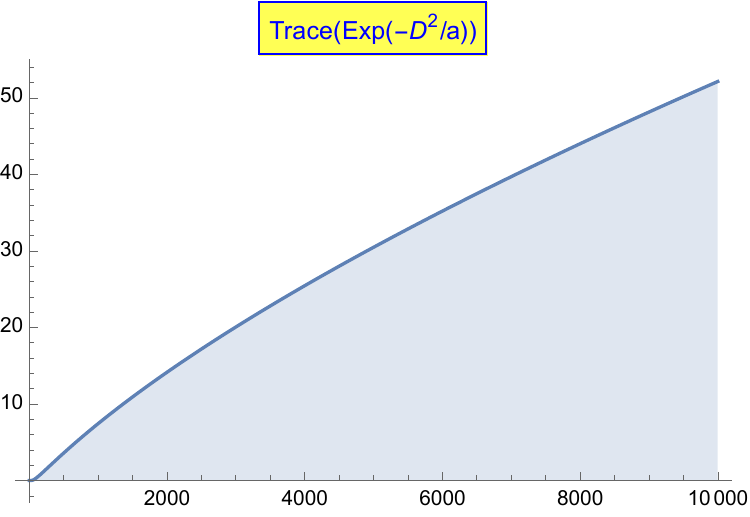}
\end{center}
\caption{Graph of $\Tr(\exp(-D^2/a))$.\label{zetaheat}}
\end{figure}

In general, the ability to compute all heat kernel coefficients in closed form is rare and typically relies on having a deep understanding of the spectral properties of the operator and the geometry of the underlying space. Most of these examples rely on spaces with high degrees of symmetry or simplicity.\newline
%The non-trivial zeros of the Riemann zeta function provide a perfect instance, assuming RH, for the inverse problem, that of finding a geometry underlying the selfadjoint operator $D$ whose spectrum is formed of the imaginary parts of non-trivial zeros of  zeta. In this paper we take a preliminary step by  computing in closed form the asymptotic expansion of $\Tr(\exp(-tD^2))$ when $t\to 0$. 
In this paper we assume $\mathrm{RH}$ and compute in closed form the asymptotic expansion of a putative self-adjoint operator whose spectrum is formed of the imaginary parts of nontrivial zeros of the Riemann zeta function as follows:
\begin{thm}\label{heatkernelintro} Assume RH and let $D$ be the self-adjoint operator whose spectrum is formed of the imaginary parts of non-trivial zeros of the Riemann zeta function. One then has the asymptotic expansion for $t\to 0$
	\begin{equation}\label{asymp}
	\Tr(\exp(-tD^2))\sim \frac{\log \left(\frac{1}{ t}\right) }{4 \sqrt{\pi } \sqrt{t}}-\frac{(\log 4\pi +\frac 12\gamma)}{2\sqrt{ \pi } \sqrt{t}}+2\exp(t/4)+\sum a_n t^{n/2}
\end{equation}
where $a_0=-\frac 14$ and for $k>0$,  using Bernouilli numbers $B_j$ and Euler numbers $E(k)$,
$$a_{2k-1}=
  \frac{\Gamma(k) \left(2^{2k-1}-1\right) B_{2 k}}{2\sqrt \pi(2 k)!}, \ \ a_{2k}=-\frac 14 \, \Gamma(k+\frac 12)\,\frac{E(2k)}{\sqrt \pi(2k)!}.
  $$
\end{thm}
The Euler numbers are defined as 
\begin{equation}\label{euler}
E(2n):=\sum _{k=1}^{2 n} \left(-\frac{1}{2}\right)^k \sum _{j=0}^{2 k} (-1)^j \binom{2 k}{j} (k-j)^{2 n}
 \end{equation}
 One has the asymptotic behavior when $k\to \infty$
$$
\frac{E(2k)}{(2k)!}\sim (-1)^k 2^{2k}\frac 4 \pi \pi^{-2k}
$$
which shows that the asymptotic expansion \eqref{asymp} is by no means convergent since its general coefficient $a_n$ diverges like a factorial.\newline
Figure \ref{zetaheat} shows the graph of $\Tr(\exp(-D^2/a))$ as a function of $a\in \R_+$.
The proof of Theorem \ref{heatkernelintro} is based on the Riemann-Weil explicit formulas which we recall in section \ref{explicit}. It is divided in two parts, we handle the contribution of the archimedean place in section \ref{archimedean} and the contribution of the primes in section \ref{primes}.\newline
The first two terms of \eqref{asymp} are the only divergent ones and their form is already a great puzzle at the level of the geometry thus revealed at the spectral level. In our joint work with H. Moscovici \cite{CM}, Theorem 5.1, we found a Dirac type square root $D$ of the prolate wave operator whose counting function for the imaginary eigenvalues  behaves as follows
 \begin{equation}\label{cm1}
 	N(E)=\frac {E}{2\pi}\left(\log \frac {E}{2\pi}-1 \right)+O(\log E).
 \end{equation}
In section \ref{counting}, we describe the link between the heat kernel asymptotic expansion and the counting function $N(E):=\#(\Sp D\cap [0,E])$ and show in Theorem \ref{estimgen} that the knowledge of the behavior of $N(E)$ in the form \eqref{cm1} suffices to secure the first terms of the heat kernel asymptotic expansion in  the form of \eqref{asymp}. This indicates that Theorem 5.1 of  \cite{CM} gives a first approximation for the geometry underlying zeta zeros.
\section{Explicit formula}\label{explicit}
Following \cite{EB}, let $f(x)$ be a smooth function on  $\mathbb{R}_{+}^*$, such that there is $\delta>0$ such that $f(x)=O\left(x^\delta\right)$ as $x \rightarrow 0+$ and $f(x)=O\left(x^{-1-\delta}\right)$ as $x \rightarrow+\infty$.
Let $\widetilde{f}(s)$ be the Mellin transform
\begin{equation}\label{mellin}
 \tilde f(s):=\int_0^\infty f(x)x^{s-1}dx.
 \end{equation}
which is an analytic function of $s$ for $-\delta<\Re(s)<1+\delta$.
 Then, with $f^\sharp(x):=x^{-1}f(x^{-1})$  the explicit formula takes the form 
 \begin{equation}\label{bombieriexplicit0}
 \sum_{\rho}\tilde f(\rho)=\int_0^\infty f(x)dx+\int_0^\infty f^\sharp(x)dx-\sum_v {\mathcal W}_v(f),
 \end{equation}
 where $v$ runs over all places  $\{\R,2,3,5,\ldots \}$ of $\Q$,  the  sum on the left hand side is over all complex zeros $\rho$ of the Riemann zeta function, and for $v=p$ 
 \begin{equation}\label{bombieriexplicit1}
 {\mathcal W}_p(f)=(\log p)\sum_{m=1}^\infty\left(f(p^m)+f^\sharp(p^m)\right).
 \end{equation}
 The archimedean distribution is defined as
 \begin{equation}\label{bombieriexplicit2}
 {\mathcal W}_\R(f):=(\log 4\pi +\gamma)f(1)+\int_{1}^\infty\left(f(x)+f^\sharp(x)-\frac 2x f(1)\right)\frac{dx}{x-x^{-1}}.
 \end{equation}
 One then has
  \begin{equation}\label{bombieriexplicit3}
 {\mathcal W}_\R(f)=(\log \pi)f(1)-\frac{1}{2\pi i}\int_{1/2+iw}\Re\left(\frac{\Gamma'}{\Gamma}\left(\frac w2\right)\right)\tilde f(w)dw.
\end{equation}
  With the above notations, let $F(x):=x^{1/2}f(x)$, one then has 
   $$F(x^{-1})=x^{-1/2}f(x^{-1})=x^{1/2}f^\sharp(x)$$
and the archimedean contribution \eqref{bombieriexplicit3} gives, using $d^*x:=dx/x$,
\begin{equation}\label{bombieriexplicit}
 W_\R(F)=(\log 4\pi +\gamma)F(1)+\int_{1}^\infty\left(F(x)+F(x^{-1})-2x^{-1/2} F(1)\right)\frac{x^{1/2}}{x-x^{-1}}d^*x.
 \end{equation}
 We consider the  Fourier transform for the duality between $\R_+^*$ and $\R$ in the form
 \begin{equation}\label{fourierm}
 {\widehat F}(s):=\int F(u)u^{-is}d^*u
 \end{equation}
 so that for $F(x):=x^{1/2}f(x)$ one has, using \eqref{mellin}, ${\widehat F}(s)=\tilde f(\frac 12-is)$. Thus assuming RH and letting $Z:=\{\rho\in \R\mid \zeta(\frac 12+i\rho)=0\}$, we can rewrite \eqref{bombieriexplicit0} as 
 \begin{equation}\label{explicit0}
 \sum_Z {\widehat F}(\rho)= {\widehat F}(i/2)+{\widehat F}(-i/2)- W_\R(F)-\sum_p W_p(F), 
 \end{equation}
 where for each prime $p$, $W_p(F)={\mathcal W}_p(f)$ which gives
 \begin{equation}\label{eachp}
 W_p(F)= (\log p)\sum_{m=1}^\infty p^{-m/2}\left(F(p^m)+F(p^{-m})\right)
 \end{equation}
\section{Archimedean contribution}\label{archimedean}

 We start from the Riemann-Weil explicit formula in the form \eqref{explicit0}
and first ignore the contributions of the finite places. With $Z:=\{\rho\in \R\mid \zeta(\frac 12+i\rho)=0\}$ as above we are thus dealing with the first three terms in the  right hand side of \eqref{explicit0}.\newline
 We choose $F_t$ such that ${\widehat F_t}(s)=\exp(-ts^2)$. One has 
 \begin{equation}\label{ft}F_t(e^y)=\frac{e^{-\frac{y^2}{4 t}}}{2\sqrt{ \pi } \sqrt{t}}, \ \ F_t(e^{-y})=F_t(e^y)
 \end{equation} Thus one gets with $ \psi(t):=\sum_p W_p(F_t)$, (see \eqref{vm} below), 
\begin{equation}\label{plicit1}
 \sum_Z \exp(-t\rho^2)= {\widehat F_t}(i/2)+{\widehat F_t}(-i/2)- W_\R(F_t)-\psi(t)=2\exp(t/4)- W_\R(F_t)-\psi(t)
 \end{equation}
  We let $x=e^u$ and get, using \eqref{bombieriexplicit},
\begin{equation}\label{plicit2}
 W_\R(F_t):=(\log 4\pi +\gamma)F_t(1)+\int_{0}^\infty\left(F_t(e^u)+F_t(e^{-u})-2e^{-u/2} F_t(1)\right)\frac{e^{u/2}}{e^u-e^{-u}}du.
 \end{equation}
 Our task now is to compute the full asymptotic expansion of \eqref{plicit2} when $t\to 0$.  One has 
 $$
  \int_0^a \left(\frac{1}{e^u-e^{-u}}-\frac{1}{2 u}\right) \, du+\int_a^{\infty } \frac{1}{e^u-e^{-u}} \, du=-\frac{1}{2} \log \left(\frac{a}{2}\right)
  $$
which taking $a=2$ gives the equality 
\begin{equation}\label{plicit3}
  \int_0^2 \left(\frac{1}{e^u-e^{-u}}-\frac{1}{2 u}\right) \, du+\int_2^{\infty } \frac{1}{e^u-e^{-u}} \, du =0
  \end{equation}
One has for $F(e^u)$ even function of $u$,
$$
I:=\int_{0}^\infty\left(F(e^u)+F(e^{-u})-2e^{-u/2} F(1)\right)\frac{e^{u/2}}{e^u-e^{-u}}du=
$$
$$
=\int_{0}^2\left(2F(e^u)-2e^{-u/2} F(1)\right)\frac{e^{u/2}}{e^u-e^{-u}}du+\int_{2}^\infty\left(2F(e^u)-2e^{-u/2} F(1)\right)\frac{e^{u/2}}{e^u-e^{-u}}du
$$
Using \eqref{plicit3} one gets
$$
\int_{2}^\infty\left(-2e^{-u/2} F(1)\right)\frac{e^{u/2}}{e^u-e^{-u}}du=\int_0^2 \left(\frac{2F(1)}{e^u-e^{-u}}-\frac{2F(1)}{2 u}\right) \, du
$$
Let 
$$I'=\int_{2}^\infty\left(2F(e^u)\right)\frac{e^{u/2}}{e^u-e^{-u}}du$$
Thus one obtains, 
$$
I=\int_{0}^2\left(\left(2F(e^u)-2e^{-u/2} F(1)\right)\frac{e^{u/2}}{e^u-e^{-u}}+\left(\frac{2F(1)}{e^u-e^{-u}}-\frac{2F(1)}{2 u}\right)\right)du+I'=
$$
$$
\int_0^2\left(2F(e^u)\frac{e^{u/2}}{e^u-e^{-u}}-\frac{2F(1)}{2 u}\right)du+I'=\int_0^2 \frac{F(e^u)-F(1)}{u}du+\int_{0}^2 2F(e^u)r(u)du+I'
$$
where
 $$r(u):=\frac{e^{u/2}}{e^u-e^{-u}}-\frac {1}{2u}$$
One has moreover
$$
I'=\int_{2}^\infty 2F(e^u) \frac{e^{u/2}}{e^u-e^{-u}}du=\int_{2}^\infty2F(e^u)\left(\frac {1}{2u}+r(u)\right)du
$$
and we obtain 
$$
I=\int_0^2 \frac{F(e^u)-F(1)}{u}du+\int_{0}^\infty 2F(e^u)r(u)du+\int_2^\infty \frac{F(e^u)}{u}du
$$
which gives the formula, for $F(e^u)$ even function of $u$,
\begin{equation}\label{wrr}
	W_\R(F)=(\log 4\pi +\gamma)F(1)+\int_0^2 \frac{F(e^u)-F(1)}{u}du+\int_{0}^\infty 2F(e^u)r(u)du+\int_2^\infty \frac{F(e^u)}{u}du
\end{equation}
We apply this formula to $F=F_t$ as in \eqref{ft}. The last term gives
$$
\int_2^\infty \frac{F_t(e^u)}{u}du=\frac{1}{2\sqrt{ \pi } \sqrt{t}}\int_2^\infty e^{-\frac{u^2}{4 t}}\frac{du}{u}=\frac{1}{2\sqrt{ \pi } \sqrt{t}}\int_{t^{-1/2}}^\infty e^{-v^2}\frac{dv}{v}=\frac{1}{4\sqrt{ \pi } \sqrt{t}}\Gamma (0,1/t)
$$
where 
 \begin{equation}\label{gam0a}
 \Gamma (0,a)=\int _a^{\infty } e^{-t}\frac{dt}{t}\sim e^{-a}/a, \ \ \text{for}\ a\to \infty
 \end{equation}
 since
 $$
 \int _a^{\infty } e^{-t}\frac{dt}{t}=e^{-a}/a \int _0^{\infty } e^{-x}\frac{dx}{1+x/a}.
 $$
 Thus when $t\to 0$, the  last term in \eqref{wrr} is $\sim \frac{\sqrt{t}}{4\sqrt{ \pi } }\exp(-1/t)=O(t^\infty)$. \newline
 We now consider the second term in \eqref{wrr}. 
 One has 
 $$
 \int_0^a \frac{1-\exp (-u)}{u} \, du=\log (a)+\Gamma (0,a)+\gamma
 $$
 
 This gives, with $u:= \frac{y^2}{4 t}$ that 
 $$
 \int_0^2 \frac{(e^{-\frac{y^2}{4 t}}-1)}{y}dy=-\frac 12\int_0^{1/t}\frac{1-\exp (-u)}{u} \, du=-\frac 12\left(\log (\frac{1}{ t})+\Gamma (0,\frac{1}{ t})+\gamma \right) 
 $$
 Thus, up to a term which is infinitely flat when $t\to 0$,
 \begin{equation}\label{as2}
\int_0^2 \frac{F_t(e^y)-F_t(1)}{y}dy=\frac{-\log \left(\frac{1}{ t}\right)-\gamma }{4 \sqrt{\pi } \sqrt{t}}+O(t^\infty)
 \end{equation}
 and the sum of the first two terms in \eqref{wrr} gives 
  \begin{equation}\label{as3}
(\log 4\pi +\gamma)F_t(1)+\int_0^2 \frac{F_t(e^u)-F_t(1)}{u}du=\frac{-\log \left(\frac{1}{ t}\right)-\gamma }{4 \sqrt{\pi } \sqrt{t}}+(\log 4\pi +\gamma)\frac{1}{2\sqrt{ \pi } \sqrt{t}}+O(t^\infty)
 \end{equation}
 It remains to deal with the third term in \eqref{wrr}, \ie $\int_{0}^\infty 2F(e^u)r(u)du$.  
 \begin{lem} One has the power series expansion,  converging for $\vert u\vert<\pi$,  \begin{equation}\label{tayll} 
 r(u)= \frac{e^{u/2}}{e^u-e^{-u}}-\frac {1}{2u}=\sum_0^\infty  b_n\, u^n,  
   \end{equation}
  where $b_0=\frac 14$ and $$b_{2k-1}=
  -\frac{ \left(1-2^{1-2k}\right) B_{2 k}}{2(2 k)!},
 \  \
  b_{2k}=\frac 14 \, 2^{-2k}\,\frac{E(2k)}{(2k)!}
  $$
 	in terms of the Bernoulli numbers $B_n$ and the Euler numbers of \eqref{euler}.
 \end{lem}
\proof  
 One has 
 $$
 \frac{e^{u/2}}{e^u-e^{-u}}=\frac{e^{u/2}}{(e^{u/2}-e^{-u/2})(e^{u/2}+e^{-u/2})}=\frac 12\left(\frac{1}{(e^{u/2}+e^{-u/2})} +\frac{1}{(e^{u/2}-e^{-u/2})} \right)
 $$
One has the Taylor expansion
$$
\frac{2}{e^{u/2}+e^{-u/2}}=1+\sum_1^\infty \frac{E(2n)}{(2n)!}\left(\frac u2\right)^{2n}
$$
 using Euler numbers defined in \eqref{euler}. 
Using Bernoulli numbers one has 
$$
\frac{2}{e^{u/2}-e^{-u/2}}=\frac 2u -\sum_1^\infty\frac{2 \left(2^{2 k-1}-1\right) B_{2 k}}{(2 k)!}\left(\frac u2\right)^{2k-1}
$$ 
One has 
$$
 B_{2 n} \sim(-1)^{n+1} 4 \sqrt{\pi n}\left(\frac{n}{\pi e}\right)^{2 n}, \ \ 
 B_{2 n}/(2n)!\sim (-1)^{n+1}  2\left(2\pi\right)^{-2n}
 $$
 and
 $$
 E_{2 n} \sim(-1)^n 8 \sqrt{\frac{n}{\pi}}\left(\frac{4 n}{\pi e}\right)^{2 n}, \  \ 2^{-2k}\,\frac{E(2k)}{(2k)!}\sim (-1)^k \frac 4 \pi \pi^{-2k}
 $$
 which checks the convergence of the series \eqref{tayll} for $\vert u\vert<\pi$  as expected from the first poles of $r(u)$ at $u=\pm i\pi$.\endproof 
\begin{lem}\label{asymplem} One has the asymptotic expansion for $t\to 0$
	 \begin{equation}\label{asymp1}
	  \int_{0}^\infty 2F_t(e^u)r(u)du\simeq	-\sum a_n t^{n/2}, \ \ a_n=-\frac{2^n \Gamma \left(\frac{n+1}{2}\right)}{\sqrt \pi}b_n
	 \end{equation}
\end{lem}
\proof We use the equality
   \begin{equation}\label{as3}
\int_0^\infty 2F_t(e^y)y^n dy=\frac{2^{n} t^{n/2} \Gamma \left(\frac{n+1}{2}\right)}{\sqrt{\pi }}
 \end{equation}
 which  gives the terms of the expansion \eqref{asymp1}.\newline
  In order to show that we get an asymptotic expansion we need to estimate the remainder in the Taylor expansion \eqref{tayll}. To avoid negative powers of $u$ we let 
  $$R(u):=2u\, r(u)=\frac{u \exp \left(\frac{u}{2}\right)}{\sinh (u)}-1$$
  and we use the Taylor formula with integral remainder for $R(u)$ at $u=0$ which gives
  \begin{equation}\label{tayll1}
 r(u)= \frac{e^{u/2}}{e^u-e^{-u}}-\frac {1}{2u}=\sum_0^k  b_n\, u^n+ \frac{1}{2u\, k!}\int_0^u R^{(k+1)}(v)(u-v)^k dv
   \end{equation}
   By construction the function $R(u)$ is smooth and moreover its derivatives fulfill
   $$
   R^{(n)}(u)\sim (-2)^{-(n-1)} (2 n-u)e^{-u/2}
   $$
   and are thus bounded on $[0,\infty)$, $\vert R^{(n)}(v)\vert\leq c_n$ which gives the inequality 
   $$
\left  \vert  \int_0^u R^{(k+1)}(v)(u-v)^k dv \right\vert \leq c_{k+1}\frac{u^{k+1}}{k+1}
   $$
   Thus the third term in \eqref{wrr} gives 
   $$
   \int_{0}^\infty 2F_t(e^u)r(u)du=\sum_0^k  b_n \int_0^\infty 2F_t(e^y)y^n dy+\rho_k(t)
   $$
   where 
   $$
   \rho_k(t)=\frac{1}{ k!}\int_{0}^\infty F_t(e^u)\left(\frac{1}{u}\int_0^u R^{(k+1)}(v)(u-v)^k dv\right)du
   $$
   so that 
   $$
   \vert \rho_k(t)\vert \leq \frac{c_{k+1}}{ (k+1)!}\int_0^\infty F_t(e^u)u^k dy=\frac{c_{k+1}\Gamma \left(\frac{k+1}{2}\right)}{\sqrt{\pi } (k+1)!}2^{k-1}t^{k/2}
   $$
   which gives the required bound for the remainder.
   \endproof 
 \section{Contribution of the primes}\label{primes}
\begin{figure}[H]	\begin{center}
\includegraphics[scale=1]{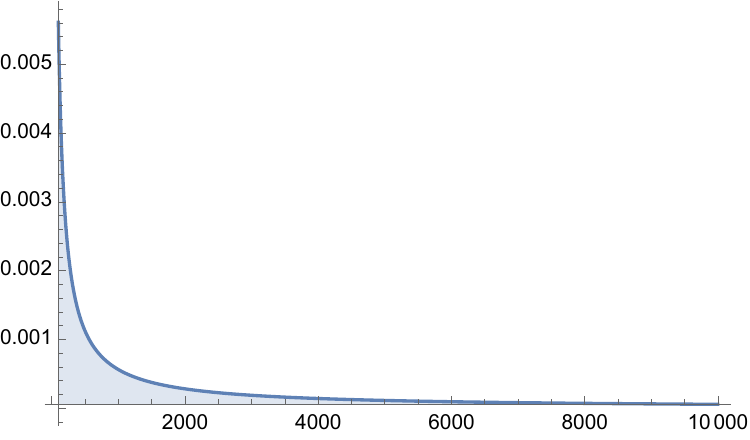}
\end{center}
\caption{Graph of discrepancy.\label{zetaheat1}}
\end{figure}

 We also need to control the contribution of the finite places and it is given, using the von Mangoldt function $\Lambda(n)$ as the  sum
 \begin{equation}\label{vm}
 \psi(t)=2\sum_2^\infty  \Lambda(n) n^{-1/2}F_t(n)=\sum_2^\infty  \Lambda(n) n^{-1/2}\,\frac{e^{-\frac{(\log n)^2}{4 t}}}{\sqrt{ \pi } \sqrt{t}}
  \end{equation}
 One has for any integer $n\geq 2$ and $t\leq t_0= \frac{\log 6}{8}\sim 0.22397$, the inequality
 $$
 e^{-\frac{(\log n)^2}{4 t}}=n^{-\frac{\log n}{4 t}}\leq 4 n^{-2} e^{-\frac{(\log 2)^2}{4 t}}
 $$
 and one thus obtains, since $\Lambda(n) n^{-1/2}\leq 1$ the  estimate for $t\leq t_0$
 $$
 \vert \psi(t)\vert \leq 4 \left(\frac{\pi ^2}{6}-1\right)\frac{e^{-\frac{(\log 2)^2}{4 t}}}{\sqrt{ \pi } \sqrt{t}}
 $$ 
 which shows that it does not contribute to the asymptotic expansion when $t\to 0$. \endproof 
 One can compute the discrepancy between $\Tr(\exp(-tD^2))$ and the first terms of the approximation \eqref{asymp} (using up to the linear term in $t$), \ie the terms
 $$
 \frac{\sqrt{\frac{1}{t}} \log \left(\frac{1}{t}\right)}{4 \sqrt{\pi }}-\frac{\gamma  \sqrt{\frac{1}{t}}}{4 \sqrt{\pi }}-\frac{\sqrt{\frac{1}{t}} \log (4 \pi )}{2 \sqrt{\pi }}+\frac{7}{4} +\frac{\sqrt t}{24 \sqrt{\pi } }+\frac{9\,t}{16}
 $$
  and for instance for $t=10^{-4}$ the difference gives 
 $-2.5\times  10^{-9}$.\newline The graph of the difference in terms of $a:=1/t\in [100,10000]$ is plotted in Figure \ref{zetaheat1} 
\section{Counting function}\label{counting}
The 
 next theorem provides the  link between the asymptotic form of the counting function (compare with Theorem 5.1 of \cite{CM}) and  the first terms of the heat expansion.  \begin{thm}\label{estimgen} Let $D$ be an operator with discrete real spectrum invariant under $x\mapsto -x$ and such that the counting function $N(E):=\#(\Sp D\cap [0,E])$
  	fulfills $N(E)=\frac {E}{2\pi}\left(\log \frac {E}{2\pi}-1 \right)+O(\log E)$, then one has for $t\to 0$,
  	\begin{equation}\label{asymp1}
	\Tr(\exp(-tD^2))= \frac{\log \left(\frac{1}{ t}\right) }{4 \sqrt{\pi } \sqrt{t}}-\frac{(\log 4\pi +\frac 12\gamma)}{2\sqrt{ \pi } \sqrt{t}}+O\left(\log (\frac{1}{ t})\right)\end{equation}
	  \end{thm}
  \proof Up to a bounded term one has, using $\partial_x(e^{-tx^2})=-2 t x e^{-t x^2}$, 
  $$
  \Tr(\exp(-tD^2))\sim 2\int_1^{\infty }\partial_EN(E)\exp(-tE^2)dE\sim 4\int_1^{\infty }N(E)tE\exp(-tE^2)dE
  $$
  Let 
  $$
 n(x):=\frac{2 x \left(\log \left(\frac{x}{2 \pi }\right)-1\right)}{2 \pi }
  $$
  In order to replace $2N(E)$ by $n(E)$ we 
 estimate the remainder due to the term in $O(\log E)$. We compute (integrating by parts) the integral
  $$
  R(t)=\int_1^{\infty } \log(x) \left(2 t x e^{-t x^2}\right) \, dx=\frac{\Gamma (0,t)}{2}=\frac{1}{2} (-\log (t)-\gamma )+O(t)
  $$
  and this suffices to give an $O(\log \frac1t)$ bound for the remainder.\newline
 It remains to compute the integral 
  \begin{align*}
    J(t)&:=\int_1^{\infty } n(x) \left(2 t x e^{-t x^2}\right) \, dx \\
    &=\frac{e^{-t}}{2 \pi  \sqrt{t}} \left(e^t G_{2,3}^{3,0}\left(t\left|
\begin{array}{c}
 1,1 \\
 0,0,\frac{3}{2} \\
\end{array}
\right.\right)-(1+\log (2 \pi )) \left(\sqrt{\pi }\, e^t\, \text{erfc}\left(\sqrt{t}\right)+2 \sqrt{t}\right)\right)
  \end{align*}
  where the Meijer G-function in this special case is the integral
  $$
  G_{2,3}^{3,0}\left(t\left|
\begin{array}{c}
 1,1 \\
 0,0,\frac{3}{2} \\
\end{array}
\right.\right)=\frac{1}{2\pi i}\int \frac{\Gamma(s)\Gamma(s)\Gamma(s+\frac 32)}{\Gamma(1+s)\Gamma(1+s)}t^{-s}ds=\frac{1}{2\pi i}\int \frac{\Gamma(s+\frac 32)}{s^2}t^{-s}ds
  $$
  The contribution of this function to the integral comes from the term in $x\log x$ in $n(x)$
  \begin{equation}\label{inti1}
  I(t):=\frac 1 \pi \int_1^{\infty } x\log x \left(2 t x e^{-t x^2}\right) \, dx
  \end{equation}
  One obtains the asymptotic expansion of $J(t)$ when $t\to 0$ in the form
  $$
 J(t)= \frac{-\log (t)-2-2 \log (2 \pi )+\Gamma'/\Gamma\left(\frac{3}{2}\right)}{4 \sqrt{\pi } \sqrt{t}}+O(t)
  $$
  To evaluate $\Gamma'/\Gamma\left(\frac{3}{2}\right)$ we use the Legendre duplication formula
  $$
  \Gamma(z)\Gamma(z+\frac 12)=2^{1-2z}\sqrt{\pi }\, \Gamma(2z)
  $$
  which gives 
  $$
  \Gamma'/\Gamma(z)+\Gamma'/\Gamma(z+\frac 12)=-2\log 2+2\,\Gamma'/\Gamma(2z)
  $$
  and thus, taking $z=1$, and using $\Gamma'/\Gamma(2)=1-\gamma$, $\Gamma'/\Gamma(1)=-\gamma$
 \begin{equation}\label{gamma}
  \Gamma'/\Gamma\left(\frac 32\right)=-2\log 2+2\,\Gamma'/\Gamma(2)-\Gamma'/\Gamma(1)=-2\log 2+2-\gamma
 \end{equation}
  We thus get   
  $$
  J(t)= \frac{-\log (t)-2-2 \log (2 \pi )+(-2\log 2+2-\gamma)}{4 \sqrt{\pi } \sqrt{t}}+O(t)=$$ $$ =\frac{\log \left(\frac{1}{ t}\right) }{4 \sqrt{\pi } \sqrt{t}}-\frac{(\log 4\pi +\frac 12\gamma)}{2\sqrt{ \pi } \sqrt{t}}+O(t)
  $$
   \endproof 
  Let us give a direct way to compute the integral \eqref{inti1} up to an $O(t)$ term for $t\to 0$. We consider 
  $$
  \int_0^\infty x^z e^{-t x^2}dx=\frac 12\,t^{-\frac z2-\frac 12}\Gamma(\frac{1+z}{2})
  $$
  and we take the derivative $\partial_z$ of both sides at $z=2$ to obtain, using $\Gamma(\frac 32)=\frac{\sqrt{\pi }}{2}$ and \eqref{gamma},
  $$
  \int_0^\infty x^2 \log(x) e^{-t x^2}dx=\frac{\sqrt{\pi }( -2\log 2+2-\gamma)}{8 t^{3/2}}-\frac{\sqrt{\pi } \log (t)}{8 t^{3/2}}
  $$
 Multiplying the result by $\frac {2t}{\pi}$ one obtains 
 $$
 I(t)=\frac{-2\log 2+2-\gamma-\log (t)}{4 \sqrt{\pi } \sqrt{t}}+O(t)
 $$ 
 The remaining term to obtain $J(t)$ is, up to a term in $O(t)$ 
 $$
 -\frac{(1+\log (2 \pi )}{\pi}\int_0^{\infty } x \left(2 t x e^{-t x^2}\right) \, dx=-\frac{1+\log (2 \pi )}{2 \sqrt{\pi } \sqrt{t}}
 $$
 which gives
 $$
J(t) =\frac{\log \left(\frac{1}{ t}\right) }{4 \sqrt{\pi } \sqrt{t}}-\frac{(\log 4\pi +\frac 12\gamma)}{2\sqrt{ \pi } \sqrt{t}}+O(t).
 $$

\end{document}